\documentclass[11pt]{article}
\usepackage[utf8]{inputenc}
\usepackage[T1]{fontenc}
\usepackage[english]{babel}
\usepackage{amsmath,amssymb,amsthm}
\usepackage{geometry}
\newtheorem{theorem}{Theorem}[section]
\newtheorem{lemma}[theorem]{Lemma}
\newtheorem{proposition}[theorem]{Proposition}
\newtheorem{corollary}[theorem]{Corollary}
\theoremstyle{definition}

\theoremstyle{remark}
\newtheorem{remark}[theorem]{Remark}
\newtheorem{example}[theorem]{Example}

\title{Potentials of measures in strictly pseudoconvex domains}
\author{Hichame Amal\footnote{Department of Mathematics, Laboratory LaREAMI, Regional Center of Trades of Education and Training, Kenitra, Morocco, hichameamal@hotmail.com}\; , Fatima Zahra Assila\footnote{Higher School of Education and Training (ESEFJ), Laboratory of Computer Science, Mathematics and Their Applications (LIMA), Science and Technology Team (ESTE),   Chouaib Doukkali University, El Jadida, Morocco,  assila.f@ucd.ac.ma}\,  and Ayoub El-Gasmi\footnote{Ibn Tofail University, Faculty of Sciences, Department of Mathematics, PO 242, Kenitra, Morocco, ayoub.el-gasmi@uit.ac.ma}}
\date{}

\begin{document}
\maketitle

\begin{abstract}
Let \(\Omega\subset\mathbb C^n\), \(n\ge2\), be a bounded strictly pseudoconvex domain.
We first show that if \(\mu\) is a finite positive measure on \(\Omega\) with no
atoms, then \((dd^c p_\mu)^n\) is absolutely continuous with respect to
Lebesgue measure. We then prove that the operation
\(u\mapsto(dd^c u)^n\mapsto p_{(dd^c u)^n}\) stabilizes the Cegrell classes
\(\mathcal E_0\) and \(\mathcal F_1\): if \(u\in\mathcal E_0(\Omega)\)
(resp.\ \(u\in\mathcal F_1(\Omega)\)), then the Green potential \(p_\mu\) of
\(\mu=(dd^c u)^n\) belongs to \(\mathcal E_0(\Omega)\) (resp.\ \(\mathcal F_1(\Omega)\)).
The \(\mathcal E_0\) statement extends a result of Carlehed from symmetric
Green functions to arbitrary strictly pseudoconvex domains. We also construct
an example with infinite-mass source whose potential lies in \(\mathcal F(\mathbb B)\),
and we show that, for any \(u\in\mathcal E_1(\mathbb B)\), the Green potential \(p_\mu\)
belongs to \(\mathcal N(\mathbb B)\) and satisfies a pointwise bound in terms of \(u\).\\

\end{abstract}

\section{Introduction}
The complex Monge--Amp\`ere operator \(u\mapsto(dd^c u)^n\) plays a central
role in K\"ahler geometry, where it encodes the equations of Calabi--Yau,
K\"ahler--Einstein and constant scalar curvature metrics, as well as in
pluripotential theory, where it is the complex analogue of the Laplace
operator. In contrast with the real theory, the operator is nonlinear for
\(n\ge2\) and discontinuous for the usual weak topologies. Its extension
beyond smooth functions was initiated by Bedford and Taylor \cite{BT76} and
developed by Cegrell \cite{Ceg98,Ceg04}.

A natural class of examples is given by the pluricomplex potentials of
positive measures,
\(
p_\mu(z)=\int_\Omega g(z,w)\,d\mu(w),
\)
where \(g\) is the pluricomplex Green function of \(\Omega\) with pole at \(w\)
and \(\mu\) is a positive finite measure. These objects, introduced
systematically by Carlehed \cite{Car99}, are the pluricomplex analogues of
the Newtonian and logarithmic potentials of classical potential theory. Among
other results, Carlehed proved that if \(p_\mu\) is bounded and
\(g(\cdot,v)\in C^2(\Omega\setminus\{v\})\) for all \(v\in\Omega\), then
\((dd^c p_\mu)^n\) is absolutely continuous with respect to Lebesgue measure,
and he established the representation formula
\[
(dd^c p_\mu)^n
=
\int_{\Omega^n}
dd^c_z g(z,v_1)\wedge\cdots\wedge dd^c_z g(z,v_n)\,
d\mu(v_1)\cdots d\mu(v_n).
\tag{1.1}
\]
For \(n\ge2\), the operation
\(\mu\mapsto p_\mu\mapsto(dd^c p_\mu)^n\) is therefore smoothing: even a
singular source \(\mu\) is mapped to an absolutely continuous Monge--Amp\`ere
measure. This phenomenon is specific to complex dimension \(n\ge2\); in
dimension \(n=1\), the operator reduces to \(dd^c p_\mu=2\pi\mu\), and no
smoothing occurs.

The singularities of pluricomplex potentials are governed by Lelong numbers.
Wiklund \cite{Wik05} showed that the Lelong number of \(p_\mu\) at \(a\)
coincides with the atom mass \(\mu(\{a\})\), and consequently that an atomless
\(\mu\) produces an atomless \((dd^c p_\mu)^n\). In general, the
Lelong--Monge--Amp\`ere inequality
\((2\pi\nu(u,x))^n\le(dd^c u)^n(\{x\})\) can be strict, pluricomplex potentials
realize the equality case, which is a distinguishing feature of this class.

A key regularity result of B\l ocki \cite{Blo00} states that, in a
\(C^\infty\) strictly pseudoconvex domain, the pluricomplex Green function
satisfies \(g(\cdot,\zeta)\in C^{1,1}(\overline\Omega\setminus\{\zeta\})\), with
second derivatives bounded up to the boundary. This regularity is optimal:
Bedford and Demailly \cite{BD88} constructed strictly pseudoconvex non-convex
domains where \(g\) fails to be \(C^2\) up to the boundary. For convex
domains, Lempert \cite{Lem81} proved that \(g\) is \(C^\infty\), but this is
a special feature of convexity. The \(C^{1,1}\) regularity is precisely the
setting where the representation formula (1.1) and the smoothing phenomenon
can be analyzed most cleanly.

The purpose of this paper is to prove a stability theorem for the operation
\(u\mapsto(dd^c u)^n\mapsto p_{(dd^c u)^n}\), and to draw some consequences
for the smoothing of Monge--Amp\`ere measures. Our main result is the
following.

\medskip
\noindent\textbf{Main Theorem.}
\emph{Let \(\Omega\subset\mathbb C^n\) be a bounded strictly pseudoconvex domain with
smooth boundary. The operation \(u\mapsto(dd^c u)^n\mapsto p_{(dd^c u)^n}\)
stabilizes the Cegrell classes \(\mathcal E_0\) and \(\mathcal F_1\): if
\(u\in\mathcal E_0(\Omega)\), then \(p_{(dd^c u)^n}\in\mathcal E_0(\Omega)\); if
\(u\in\mathcal F_1(\Omega)\), then \(p_{(dd^c u)^n}\in\mathcal F_1(\Omega)\). In both cases,
\((dd^c p_{(dd^c u)^n})^n\) is absolutely continuous with respect to
Lebesgue measure.}
\medskip

The \(\mathcal E_0\) statement extends Carlehed's Proposition~9.3, which
requires the pluricomplex Green function to be symmetric, to arbitrary bounded
strictly pseudoconvex domains with \(C^\infty\) boundary. In \cite{Ceg12}, Cegrell proved that \(p_\nu\in\mathcal F_1(\Omega)\) under a
different hypothesis on the source measure (see Remark~\ref{rem:R1}). The absolute
continuity statement is a consequence of the smoothing theorem proved in
Section~\ref{sec:mixed}, which itself extends Carlehed's Theorem~8.1 to
strictly pseudoconvex domains without any boundedness assumption on
\(p_\mu\). We also construct, in Section~\ref{sec:infinite-mass}, an explicit
example where the source measure has infinite total mass while its Green
potential belongs to \(\mathcal F(\mathbb B)\), showing that the smoothing phenomenon
persists in the infinite-mass regime.

The paper is organized as follows. Section~\ref{sec:prelim} collects
preliminaries on plurisubharmonic functions, Cegrell classes, capacities, and
the pluricomplex Green function. Section~\ref{sec:mixed} is devoted to mixed
measures and atoms: we prove that
\(dd^c g(z,v_1)\wedge\cdots\wedge dd^c g(z,v_n)\) is absolutely continuous
whenever the poles are not all equal, compute the mass of the atoms of
\((dd^c p_\mu)^n\), and deduce the smoothing theorem
(Theorem~\ref{thm:bdd}). Section~\ref{sec:E1} contains the stability theorem
(Theorem~\ref{thm:stability}) and its two proofs.
Section~\ref{sec:infinite-mass} treats the infinite-mass case: it contains an
explicit construction with \(\mu(\mathbb B)=+\infty\) and \(p_\mu\in\mathcal F(\mathbb B)\), and a
general result stating that \(p_\mu\in\mathcal N(\mathbb B)\) for every \(u\in\mathcal E_1(\mathbb B)\),
together with a pointwise bound.

\section{Preliminaries}
\label{sec:prelim}
Let \(\Omega\subset\mathbb C^n\) be an open domain. We set
\[
d^c:=\frac{i}{2\pi}\bigl(\bar\partial-\partial\bigr),
\qquad
dd^c:=\frac{i}{\pi}\,\partial\bar\partial,
\qquad
\beta:=dd^c|z|^2.
\]

\noindent{\bf Plurisubharmonic functions.}
A function \(u:\Omega\to[-\infty,+\infty[\) is \emph{plurisubharmonic} (psh) if
it is upper semicontinuous, not identically \(-\infty\), and restricts to a
subharmonic function on every complex line; equivalently, \(dd^c u\ge0\) in
the sense of positive currents. We denote by \(\operatorname{PSH}(\Omega)\) the convex cone of
psh functions and by \(\operatorname{PSH}^-(\Omega)\) the negative ones. The cone
\(\operatorname{PSH}(\Omega)\) is contained in \(L^1_{\mathrm{loc}}(\Omega)\), stable under
finite maxima and decreasing limits, and closed under \(L^1_{\mathrm{loc}}\)
convergence of locally upper bounded sequences (Hartogs' lemma, Montel
property).

A set \(P\subset\Omega\) is \emph{pluripolar} if it is locally contained in the
polar set \(\{u=-\infty\}\) of some \(u\in\operatorname{PSH}(\Omega)\); pluripolar sets have
zero Lebesgue measure and do not charge any Monge--Amp\`ere measure
\cite{Kli91}. A domain \(\Omega\) is \emph{hyperconvex} if it admits a negative
psh exhaustion function, and \emph{strictly pseudoconvex} if this exhaustion
can be taken in \(\operatorname{PSH}^-(\Omega)\cap C^\infty(\overline\Omega)\) with
non-vanishing gradient on \(\partial\Omega\).
\\

\noindent{\bf The complex Monge--Amp\`ere operator.}
For \(u\in\operatorname{PSH}(\Omega)\cap L^\infty_{\mathrm{loc}}(\Omega)\), the complex
Monge--Amp\`ere operator \(u\mapsto(dd^c u)^n\) is defined by the
Bedford--Taylor theory \cite{BT76}: for smooth \(u\),
\((dd^c u)^n=c_n\det(\partial^2 u/\partial z_j\partial\bar z_k)\,dV\), and
one extends recursively by
\(dd^c u_1\wedge\cdots\wedge dd^c u_k
:=dd^c(u_1\,dd^c u_2\wedge\cdots\wedge dd^c u_k)\); the result is a positive
Radon measure independent of the order of \(u_1,\dots,u_k\). The operator is
continuous along monotone sequences of locally bounded psh functions, but not
for the \(L^1_{\mathrm{loc}}\)-topology. We shall use the
Chern--Levine--Nirenberg estimate \cite{Kli91}: for every compact
\(K\Subset\Omega\) there exist \(L\Subset\Omega\) and \(C>0\) such that
\(\int_K dd^c u_1\wedge\cdots\wedge dd^c u_n\le C\prod_i\|u_i\|_{L^\infty(L)}\).
\\

\noindent{\bf Cegrell classes.}
Following Cegrell \cite{Ceg98,Ceg04}, define
\[
\mathcal E_0(\Omega)
:=
\left\{u\in\operatorname{PSH}(\Omega)\cap L^\infty(\Omega):\lim_{z\to\partial\Omega}u(z)=0,\ 
\int_\Omega(dd^c u)^n<\infty\right\}.
\]
The class \(\mathcal E(\Omega)\subset\operatorname{PSH}^-(\Omega)\) consists of all \(u\) such that, for
every \(z_0\in\Omega\), there exist a neighborhood \(\omega\) of \(z_0\) and a
decreasing sequence \((h_j)\subset\mathcal E_0(\Omega)\) with \(h_j\searrow u\)
on \(\omega\) and \(\sup_j\int_\Omega(dd^c h_j)^n<\infty\). Finally,
\[
\mathcal F(\Omega):=\left\{u\in\mathcal E(\Omega):
\begin{array}{l}
\exists\,(u_j)\subset\mathcal E_0(\Omega),\ u_j\searrow u\ \text{on }\Omega,\\[2pt]
\sup_j\int_\Omega(dd^c u_j)^n<+\infty
\end{array}\right\}.
\]
For \(p>0\), \(\mathcal E_p(\Omega):=\{u\in\mathcal E(\Omega):\int_\Omega|u|^p(dd^c u)^n<\infty\}\).
We shall use the class
\[
\mathcal F_1(\Omega):=\mathcal F(\Omega)\cap\mathcal E_1(\Omega),
\]
of functions in \(\mathcal F\) with finite energy \(\int_\Omega-u\,(dd^c u)^n\) and
finite Monge--Amp\`ere mass.

The class \(\mathcal E(\Omega)\) is the largest class of negative psh functions on
which the operator is well defined and continuous under decreasing limits.
On \(\mathcal F(\Omega)\), integration by parts holds \cite[Corollary~3.4]{Ceg04}: for
\(u,v\in\mathcal F(\Omega)\) and \(T\) closed positive of bidegree \((1,1)\),
\(\int_\Omega v\,dd^c u\wedge T=\int_\Omega u\,dd^c v\wedge T\), provided one of
the two integrals is finite. We shall also use the monotonicity property:
if \(u\le v\) both belong to \(\mathcal F(\Omega)\) and \(\phi\in\operatorname{PSH}^-(\Omega)\), then
\[
\int_\Omega-\phi\,(dd^c v)^n\le\int_\Omega-\phi\,(dd^c u)^n.
\tag{2.1}
\]
The comparison principle holds: if
\(u,v\in\operatorname{PSH}(\Omega)\cap L^\infty(\Omega)\), \(\liminf_{z\to\partial\Omega}(u-v)\ge0\),
and \((dd^c u)^n\le(dd^c v)^n\), then \(u\ge v\) \cite[Theorem~5.12]{Ceg04}.
Finally, we recall the \emph{localization maximum principle} of Khue--Hiep
\cite[Theorem~4.1]{KHH09}: for
\(u,u_1,\dots,u_{n-1}\in\mathcal E(\Omega)\), \(v\in\operatorname{PSH}^-(\Omega)\), and
\(T:=dd^c u_1\wedge\cdots\wedge dd^c u_{n-1}\), one has
\[
dd^c\max(u,v)\wedge T|_{\{u>v\}}
=
dd^c u\wedge T|_{\{u>v\}}.
\tag{2.2}
\]
In the particular case \(u_1=\cdots=u_{n-1}=u\), iterating (2.2) on the
open set \(\{u>v\}\) where \(\max(u,v)=u\) gives the classical identity
\(
\mathbf{1}_{\{u>v\}}(dd^c\max(u,v))^n
=
\mathbf{1}_{\{u>v\}}(dd^c u)^n.
\)
\\

\noindent{\bf The class \(\mathcal N(\Omega)\).}
Let \((\Omega_j)\) be a fundamental sequence of \(\Omega\) (increasing, strictly
pseudoconvex, relatively compact, with union \(\Omega\)). For \(u\in\mathcal E(\Omega)\),
put \(u^j:=\sup\{\phi\in\operatorname{PSH}(\Omega):\phi\le u\ \text{on }\Omega\setminus\Omega_j\}\)
and \(\tilde u:=(\lim_j u^j)^*\). The function \(\tilde u\) does not depend on
\((\Omega_j)\), belongs to \(\mathcal E(\Omega)\), is maximal (\((dd^c\tilde u)^n=0\)
\cite{Blo04,Blo06}), and satisfies \(u\le\tilde u\le0\). Define
\[
\mathcal N(\Omega):=\{u\in\mathcal E(\Omega):\tilde u=0\}.
\]
One has \(\mathcal E_0\subset\mathcal F\subset\mathcal N\subset\mathcal E\). The characterization of
Czy\.z--Hed \cite[Proposition~2.2(c)]{CzyzHed} will be crucial below:
\[
u\in\mathcal N(\Omega)\ \text{and}\ \int_\Omega(dd^c u)^n<+\infty
\quad\Longrightarrow\quad
u\in\mathcal F(\Omega).
\tag{2.3}
\]

\noindent{\bf Capacity and Cegrell's H\"older inequality.}
The \emph{Monge--Amp\`ere capacity} of a Borel set \(E\subset\Omega\) is
\[
\operatorname{Cap}_\Omega(E)
:=
\sup\left\{\int_E(dd^c u)^n:\ u\in\operatorname{PSH}(\Omega),\ -1\le u\le0\right\}.
\]
Pluripolar sets are exactly those of zero outer capacity \cite{Kli91}, and
psh functions are quasi-continuous with respect to this capacity. We shall
use Cegrell's H\"older inequality \cite[Theorem~5.5]{Ceg04}: for
\(u_1,\dots,u_n\in\mathcal F(\Omega)\) and \(h\in\mathcal E_0(\Omega)\),
\[
\int_\Omega -h\,dd^c u_1\wedge\cdots\wedge dd^c u_n
\le
\prod_{i=1}^n
\left(\int_\Omega -h\,(dd^c u_i)^n\right)^{1/n}.
\]
\\

\noindent{\bf Lelong numbers.}
For \(u\in\operatorname{PSH}(\Omega)\) and \(a\in\Omega\), the \emph{Lelong number} of $u$ at $a$ is
\(\nu(u,a):=\lim_{r\to0^+}\sup_{B(a,r)}u/\log r\). 
For \(u\in\mathcal E(\Omega)\), the sharp estimate
\(2\pi\nu(u,a)\le((dd^c u)^n(\{a\}))^{1/n}\) holds
\cite[Corollary~5.7]{Ceg04}; in particular, \(\{\nu(u,\cdot)>0\}\) is
discrete. We shall use Demailly's comparison theorem \cite{Dem87}: if
\(u,v\in\operatorname{PSH}(\Omega)\) satisfy \(u(z)\le v(z)+O(1)\) as \(z\to a\) and
\(u(z)/v(z)\to1\), then \((dd^c u)^n(\{a\})\le(dd^c v)^n(\{a\})\), with
equality if the upper limit is a limit.
\\

\noindent{\bf The pluricomplex Green function.}
For \(a\in\Omega\), the \emph{pluricomplex Green function} of \(\Omega\) with pole
at \(a\) is
\[
g_\Omega(z,a):=\sup\left\{u(z):
\begin{array}{l}
u\in\operatorname{PSH}(\Omega),\ u<0,\\[2pt]
u(z)\le\log|z-a|+O(1)\ \text{as }z\to a
\end{array}\right\}.
\]
It is negative, maximal, has a logarithmic singularity at \(a\), and
vanishes on \(\partial\Omega\). The normalization
\((dd^c g_\Omega(\cdot,a))^n=(2\pi)^n\delta_a\) has been fixed above. For the
unit ball \(\mathbb B\subset\mathbb C^n\), \(g_\mathbb B(z,a)=\log|\Phi_a(z)|\) with \(\Phi_a\)
the automorphism of \(\mathbb B\) sending \(a\) to \(0\).

\begin{theorem}[B\l ocki \cite{Blo00}]
\label{thm:green-C11}
Let \(\Omega\subset\mathbb C^n\) be a \(C^\infty\) strictly pseudoconvex domain and
\(g=g_\Omega(\cdot,\zeta)\). Then
\(g\in C^{1,1}(\overline\Omega\setminus\{\zeta\})\).
\end{theorem}
The \(C^{1,1}\) regularity of Theorem~\ref{thm:green-C11} is optimal. Let
\(
C:=\{(x_1,x_2,0,0)\in\mathbb C^2:x_1^2+x_2^2=1\}
\)
be the unit circle in the real \(2\)-plane \(\mathbb R^2\subset\mathbb C^2\), where
we write \((z,w)=(x_1+iy_1,x_2+iy_2)\). For small \(\varepsilon>0\), set
\[
\Omega_\varepsilon
:=\{(z,w)\in\mathbb C^2:\operatorname{dist}((z,w),C)<\varepsilon\}.
\]
Then \(\Omega_\varepsilon\) is a bounded strictly pseudoconvex non-convex
domain with \(C^\infty\) boundary. Fix \(p=(1,0)\). By
Theorem~\ref{thm:green-C11},
\(g_{\Omega_\varepsilon}(\cdot,p)\in C^{1,1}\) up to the boundary, while
Bedford and Demailly \cite{BD88} proved that, for suitable \(\varepsilon\),
\(g_{\Omega_\varepsilon}(\cdot,p)\notin C^2\). Thus the \(C^{1,1}\)
regularity cannot be improved to \(C^2\).
 For convex domains, Lempert \cite{Lem81} proved
that \(g\in C^\infty\).
\\

\noindent{\bf Potentials of measures.}
For a finite positive measure \(\mu\) on \(\Omega\), the \emph{pluricomplex
potential} (or \emph{Green potential}) of \(\mu\) is
\[
p_\mu(z)=\int_\Omega g_\Omega(z,a)\,d\mu(a).
\]
It satisfies \(p_\mu\in\operatorname{PSH}(\Omega)\), \(p_\mu\le0\), and if \(\Omega\) is hyperconvex,
then \(p_\mu\) vanishes on the boundary. Moreover \(p_\mu\in\mathcal F(\Omega)\) by Cegrell's theorem
\cite[Theorem~1.1]{Ceg12} (see also Wiklund \cite[Lemma~4.3]{Wik05}). The remaining part of
\cite[Theorem~1.1]{Ceg12} gives the sharp estimate
\[
\int_\Omega(-h)\,(dd^c p_\mu)^n
\le
\left(\int_\Omega(-h)^{1/n}\,d\mu\right)^n\le \mu(\Omega)^{n-1}\int_\Omega -hd\mu
\qquad
\text{for all }h\in\operatorname{PSH}^-(\Omega).
\]
If \(n\ge2\), \(g(.,w)\in\mathcal{C}^2(\Omega\setminus\{w\})\) and \(p_\mu\) is bounded, then \((dd^c p_\mu)^n\) is absolutely
continuous with respect to Lebesgue measure \cite[Theorem~8.1]{Car99}.

\section{Mixed measures, atoms, and examples}
\label{sec:mixed}

Throughout this section, \(\Omega\subset\mathbb C^n\), \(n\ge2\), is a bounded strictly
pseudoconvex domain with smooth boundary, \(g(z,w)\) is its pluricomplex
Green function, and \(\lambda\) is the Lebesgue measure on \(\mathbb C^n\). By
Theorem~\ref{thm:green-C11}, \(g(\cdot,v)\in C^{1,1}(\Omega\setminus\{v\})\).
For a finite positive measure \(\mu\), put
\(p_\mu(z)=\int_\Omega g(z,w)\,d\mu(w)\).

\begin{proposition}
\label{prop:mixed}
Let \(v_1,\dots,v_n\in\Omega\) and
\(\nu_{v_1,\dots,v_n}:=dd^c_z g(z,v_1)\wedge\cdots\wedge dd^c_z g(z,v_n)\).
If the points \(v_1,\dots,v_n\) are not all equal, then
\(\nu_{v_1,\dots,v_n}\ll\lambda\). If \(v_1=\cdots=v_n=v\), then
\(\nu_{v,\dots,v}=(2\pi)^n\delta_v\).
\end{proposition}

\begin{proof}
Let \(A=\{v_1,\dots,v_n\}\). Outside \(A\), each \(g(\cdot,v_j)\) is
\(C^{1,1}\), so \(\nu_{v_1,\dots,v_n}=F\,d\lambda\) on \(\Omega\setminus A\) with
\(F\in L^\infty_{\mathrm{loc}}\). The singular part of \(\nu_{v_1,\dots,v_n}\)
is supported on \(A\); we show it vanishes. Fix \(a\in A\), let \(1\leq k<n\) be the
multiplicity of \(a\) among the poles and \(b_1,\dots,b_{n-k}\) the other poles. We must
show that
\[
I_a:=\int_{\{a\}}\bigl(dd^c g(z,a)\bigr)^k\wedge dd^c g(z,b_1)\wedge\cdots\wedge dd^c g(z,b_{n-k})=0.
\]
Set \(u(z):=g(z,a)+\sum_j g(z,b_j)\). Near \(a\), \(u(z)\le g(z,a)\) and
\(u(z)/g(z,a)\to1\) as \(z\to a\). Demailly's comparison theorem
\cite[Theorem~5.1]{Dem93} gives
\(\int_{\{a\}}(dd^c u)^n=\int_{\{a\}}(dd^c g(\cdot,a))^n=(2\pi)^n\).
Expanding \((dd^c u)^n\) by the multinomial formula, all coefficients being
positive, and using the positivity of the mixed Monge--Amp\`ere measures, we
obtain
\[
(dd^c u)^n\ge\bigl(dd^c g(z,a)\bigr)^n+\bigl(dd^c g(z,a)\bigr)^k\wedge dd^c g(z,b_1)\wedge\cdots\wedge dd^c g(z,b_{n-k}),
\]
so integrating over \(\{a\}\) yields \((2\pi)^n\ge(2\pi)^n+I_a\), hence
\(I_a=0\).
\end{proof}

\begin{lemma}
\label{lem:representation}
Let \(\Omega\subset\mathbb C^n\) be a bounded hyperconvex domain, \(\mu\) a positive
finite measure, and \(p_\mu(z)=\int_\Omega g(z,w)\,d\mu(w)\). Then
\[
(dd^c p_\mu)^n
=
\int_{\Omega^n}
dd^c_z g(z,v_1)\wedge\cdots\wedge dd^c_z g(z,v_n)\,
d\mu(v_1)\cdots d\mu(v_n),
\tag{3.1}
\]
as an identity of positive Radon measures.
\end{lemma}

\begin{proof}
By Cegrell's theorem \cite[Theorem~1.1]{Ceg12}, since \(\mu\) is finite,
\(p_\mu\in\mathcal F(\Omega)\) and \(\int_\Omega(dd^c p_\mu)^n\le\mu(\Omega)^n<\infty\). For
each \(v\in\Omega\), the Green function \(g(\cdot,v)\in\mathcal F(\Omega)\)
 and its Monge--Amp\`ere measure is \((2\pi)^n\delta_v\).
Thus, for every \((v_1,\dots,v_n)\), the mixed measure
\(dd^c g(\cdot,v_1)\wedge\cdots\wedge dd^c g(\cdot,v_n)\) is well defined and
positive. By Carlehed \cite[Section~6]{Car99}, its total mass is uniformly
bounded in \((v_1,\dots,v_n)\), so the right-hand side of (3.1), which we
denote by \(\nu\), is a well-defined positive Radon measure.
It suffices to prove \((dd^c p_\mu)^n=\nu\). By Cegrell's Lemma~3.1
\cite[Lemma~3.1]{Ceg04}, every \(\chi\in C_0^\infty(\Omega)\) is a difference of
two functions in \(\mathcal E_0(\Omega)\cap C(\overline\Omega)\). By linearity, it
is enough to show that
\[
\int_\Omega h\,(dd^c p_\mu)^n=\int_\Omega h\,d\nu
\tag{3.2}
\]
for every \(h\in\mathcal E_0(\Omega)\).

Fix \(h\in\mathcal E_0(\Omega)\) and \((v_1,\dots,v_n)\in\Omega^n\), and write
\(g_j:=g(\cdot,v_j)\). Since \(h,g_j\in\mathcal F(\Omega)\), integration by parts in
\(\mathcal F(\Omega)\) \cite[Corollary~3.4]{Ceg04} yields
\[
\int_\Omega h\,dd^c g_1\wedge\cdots\wedge dd^c g_n
=
\int_\Omega g_1\,dd^c h\wedge dd^c g_2\wedge\cdots\wedge dd^c g_n .
\tag{3.3}
\]
All the integrals are finite because the mixed masses are uniformly bounded
and \(h\in L^\infty\). Now integrate (3.3) against
\(d\mu(v_1)\cdots d\mu(v_n)\). By Fubini, which is justified by the
positivity and the uniform bound on the mixed masses, the left-hand side is
\(\int_\Omega h\,d\nu\); on the right-hand side, integrating first in \(v_1\)
and using \(\int_\Omega g_1\,d\mu(v_1)=p_\mu\) gives
\begin{multline*}
\int_{\Omega^n}\!\Bigl(\int_\Omega g_1\,dd^c h\wedge dd^c g_2\wedge\cdots\wedge dd^c g_n\Bigr)\,d\mu(v_1)\cdots d\mu(v_n)\\
=
\int_{\Omega^{n-1}}\!\Bigl(\int_\Omega p_\mu\,dd^c h\wedge dd^c g_2\wedge\cdots\wedge dd^c g_n\Bigr)\,d\mu(v_2)\cdots d\mu(v_n).
\end{multline*}
Applying integration by parts to \((p_\mu,g_2)\) with current
\(T_2:=dd^c h\wedge dd^c g_3\wedge\cdots\wedge dd^c g_n\) and integrating
against \(d\mu(v_2)\), we obtain
\(\int_\Omega p_\mu\,dd^c h\wedge dd^c p_\mu\wedge T_2\), back to the original
configuration with one \(g_j\) replaced by \(p_\mu\). Iterating \(n-1\)
times gives \(\int_\Omega p_\mu\,dd^c h\wedge (dd^c p_\mu)^{n-1}\), and a
final integration by parts between \(p_\mu\) and \(h\) converts this into
\(\int_\Omega h\,(dd^c p_\mu)^n\). Hence (3.2) holds for every
\(h\in\mathcal E_0(\Omega)\) and then \((dd^c p_\mu)^n=\nu\).
\end{proof}

\begin{corollary}
\label{cor:atoms}
Let \(n\ge2\) and let \(\mu\) be a positive finite measure on \(\Omega\). Then
\((dd^c p_\mu)^n\) has an atom at \(a\) if and only if \(\mu\) has an atom at
\(a\), and
\[
(dd^c p_\mu)^n(\{a\})=(2\pi)^n\mu(\{a\})^n=(2\pi\nu(p_\mu,a))^n.
\]

\end{corollary}

\begin{proof}
By Lemma~\ref{lem:representation}, (3.1) holds. Combined with
Proposition~\ref{prop:mixed}, the mixed measure
\(dd^c_z g(z,v_1)\wedge\cdots\wedge dd^c_z g(z,v_n)\) is absolutely
continuous when the poles are not all equal, and equals
\((2\pi)^n\delta_a\) when \(v_1=\cdots=v_n=a\). Integrating against
\(\mu^{\otimes n}\), the only contribution to the atom at \(a\) is the
diagonal, giving \((2\pi)^n\mu(\{a\})^n\). The Lelong formula is
\cite[Corollary~4.5]{Wik05}.
\end{proof}

\begin{theorem}
\label{thm:bdd}
Let \(n\ge2\), let \(\mu\) be a positive finite measure on \(\Omega\) with no
atoms. Then \((dd^c p_\mu)^n\ll\lambda\).
\end{theorem}

\begin{proof}
Let \(N\subset\Omega\) be \(\lambda\)-null. By Lemma~\ref{lem:representation} and
Fubini,
\[
\int_N(dd^c p_\mu)^n=\int_{\Omega^n}\left(\int_N dd^c_z g(z,v_1)\wedge\cdots\wedge dd^c_z g(z,v_n)\right)d\mu(v_1)\cdots d\mu(v_n).
\]
By Proposition~\ref{prop:mixed}, the inner integral vanishes when the poles
are not all equal, and equals \((2\pi)^n\mathbf{1}_{N}(v)\) on the diagonal.
Hence \(\int_N(dd^c p_\mu)^n\le(2\pi)^n\int_\Omega\mu(\{w\})^{n-1}d\mu(w)=0\).
\end{proof}
\begin{remark}
For \(n=1\), \(dd^c p_\mu=2\pi\mu\), so \((dd^c p_\mu)\) is absolutely
continuous iff \(\mu\) is; the smoothing phenomenon
\(\mu\mapsto p_\mu\mapsto(dd^c p_\mu)^n\) is specific to \(n\ge2\).
\end{remark}

\begin{example}[Unbounded potential from an explicit density]
\label{ex:density}
Let \(n\ge2\) and \(B=B(0,1/2)\subset\mathbb C^n\). Define the positive measure
\[
  d\mu(w) := \frac{\mathbf 1_{B}(w)}{|w|^{2n}\,(-\log|w|)^2}\,d\lambda(w).
\]
Then the following hold.

\begin{enumerate}
\item 
Writing \(w=r\xi\) with \(r=|w|>0\) and \(|\xi|=1\), the Lebesgue measure
decomposes as \(d\lambda(w)=c_n\,r^{2n-1}dr\,d\sigma(\xi)\) for a dimensional constant
\(c_n>0\), where \(d\sigma\) is the normalized surface measure on the unit sphere.
Thus
\[
  \mu(B)
  = c_n\int_0^{1/2}\frac{r^{2n-1}}{r^{2n}(-\log r)^2}\,dr
  = c_n\int_0^{1/2}\frac{1}{r(-\log r)^2}\,dr = \frac{c_n}{\log 2}<+\infty.
\]
The density is locally integrable, hence \(\mu\ll\lambda\) and \(\mu\) is atomless. Moreover \(p_\mu\in\mathcal F(\Omega)\) by Cegrell's theorem
\cite[Theorem~1.1]{Ceg12}.

\item 
For the ball \(B(0,1/2)\), the pluricomplex Green function with pole at \(0\) is
\(g_B(0,w)=\log(2|w|)\). Therefore
\[
  p_\mu(0)
  = \int_B \log(2|w|)\,d\mu(w)
  = \mu(B)\log 2 + \int_B \log|w|\,d\mu(w).
\]
The first term is finite. For the second term, using the same polar coordinates,
\[
  \int_B \log|w|\,d\mu(w)
  = c_n\int_0^{1/2}\frac{\log r}{r(-\log r)^2}\,dr
  = -c_n\int_0^{1/2}\frac{1}{r(-\log r)}\,dr
  = -c_n\int_{\log 2}^{+\infty}\frac{dt}{t}=-\infty .
\]

Thus \(p_\mu(0)=-\infty\), and \(p_\mu\) is not bounded below on \(B\).

\item Since \(\mu\) is a finite atomless positive measure on \(B\), Theorem~\ref{thm:bdd}
applies and yields \((dd^c p_\mu)^n\ll\lambda\). In particular,
\((dd^c p_\mu)^n(\{0\})=0\), even though \(p_\mu(0)=-\infty\).

\end{enumerate}
\end{example}

\section{Stability of the classes \(\mathcal E_0\) and \(\mathcal F_1\) under the Green potential}
\label{sec:E1}

We show that the operation
\[
u\longmapsto\mu:=(dd^c u)^n\longmapsto p_\mu
\]
leaves the Cegrell classes \(\mathcal E_0\) and \(\mathcal F_1\) invariant: if \(\mu\)
is the Monge--Amp\`ere measure of a function in \(\mathcal E_0\) (resp.\
\(\mathcal F_1\)), then its Green potential \(p_\mu\) also belongs to \(\mathcal E_0\)
(resp.\ \(\mathcal F_1\)). The two cases are gathered in a single statement, proved
by different techniques: the \(\mathcal E_0\) case uses a bounded psh
dominator and Cegrell's existence theorem \cite[Proposition~6.1]{Ceg98},
whereas the \(\mathcal F_1\) case uses energy estimates and the comparison principle.

\begin{theorem}
\label{thm:stability}
Let \(\Omega\subset\mathbb C^n\) be a bounded strictly pseudoconvex domain with smooth
boundary, \(u\in\mathcal E(\Omega)\)  and
\(\mu=(dd^c u)^n\). Then:
\begin{enumerate}
\item[(a)] if \(u\in\mathcal E_0(\Omega)\), then \(p_\mu\in\mathcal E_0(\Omega)\);
\item[(b)] if \(u\in\mathcal F_1(\Omega)\), then \(p_\mu\in\mathcal F_1(\Omega)\).
\end{enumerate}
In both cases, \((dd^c p_\mu)^n\) is absolutely continuous with respect to
Lebesgue measure.
\end{theorem}

The absolute continuity statement follows in both cases from
Theorem~\ref{thm:bdd}: if \(u\in\mathcal E_0(\Omega)\) or \(u\in\mathcal F_1(\Omega)\),
then \(u\in\mathcal E_1(\Omega)\), hence \(\mu=(dd^c u)^n\) does not charge pluripolar
sets and is atomless. We therefore focus on the membership statements.

\subsection{Proof of Theorem~\ref{thm:stability}(a): the \(\mathcal E_0\) case}

Let \(u\in\mathcal E_0(\Omega)\), \(\mu=(dd^c u)^n\), and
\(p:=p_\mu\). We prove that \(p\in\mathcal E_0(\Omega)\), i.e.\ that \(p\) is
bounded, vanishes on \(\partial\Omega\), and has finite Monge--Amp\`ere mass.

\medskip\noindent\textbf{Step 1: boundedness and boundary behaviour.}
Since \(g_\Omega\le0\), we have \(p\le0\). We prove a lower bound.
Let \(R>0\) with \(\Omega\Subset B(R)\subset\mathbb C^n\). Fix \(K\Subset\Omega\).
Choose \(r_2>0\) with \(K+\overline{B(0,r_2)}\Subset\Omega\), set
\(r_1:=r_2/2\), and fix \(z\in K\). Put
\[
A:=\frac{\|u\|_\infty}{\log 2},\qquad j:=\|u\|_\infty+1,
\]
and define, on \(B(R)\setminus\{z\}\),
\[
v^{(z)}_2(w):=\max\Bigl\{A\log\frac{|w-z|}{r_2},\,-j\Bigr\}.
\]
Then \(v^{(z)}_2\in\operatorname{PSH}(B(R))\cap L^\infty(B(R))\), with
\[
-j\le v^{(z)}_2(w)\le A\log\frac{2R}{r_2}
\qquad\text{for all }w\in B(R)\setminus\{z\}.
\]
Moreover, for \(|w-z|=r_2\), \(v^{(z)}_2(w)=0>u(w)\), and for
\(|w-z|\le r_1=r_2/2\),
\[
A\log\frac{|w-z|}{r_2}\le A\log\frac12=-\|u\|_\infty\le u(w).
\]
Set
\[
W_z(w):=
\begin{cases}
\max\{v^{(z)}_2(w),\,u(w)\}, & w\in B(z,r_2),\\
v^{(z)}_2(w), & w\in B(R)\setminus B(z,r_2).
\end{cases}
\]
Since \(v^{(z)}_2=0>u\) on \(\partial B(z,r_2)\), the gluing principle gives
\(W_z\in\operatorname{PSH}(B(R))\cap L^\infty(B(R))\), and \(W_z=u\) on \(B(z,r_1)\).
In particular,
\[
(dd^c W_z)^n\ge\mathbf 1_{B(z,r_1)}\,\mu\quad\text{on }B(R),
\qquad
W_z(z)=u(z).
\tag{4.1}
\]

Set
\[
\Lambda:=A\log\frac{2R}{r_2},\qquad
\widetilde W_z:=W_z-\Lambda.
\]
Since \(W_z\le\Lambda\) on \(B(R)\), we have \(\widetilde W_z\le0\), and
\(
(dd^c\widetilde W_z)^n=(dd^c W_z)^n\ge\mathbf 1_{B(z,r_1)}\,\mu.
\)
By \cite[Proposition~6.1]{Ceg98} and its proof, applied to the bounded psh function
\(\widetilde W_z\), there exists \(\varphi_z\in\mathcal F_1(B(R))\) with
\[
(dd^c\varphi_z)^n=\mathbf 1_{B(z,r_1)}\,\mu,
\qquad
\varphi_z\ge\widetilde W_z.
\tag{4.2}
\]
In particular, using \(W_z(z)=u(z)\),
\[
\varphi_z(z)\ge\widetilde W_z(z)=u(z)-\Lambda\ge-\|u\|_\infty-\Lambda.
\tag{4.3}
\]

Since \(\varphi_z\in\mathcal F_1(B(R))\subset\mathcal F(B(R))\) and
\(g_{B(R)}(z,\cdot)\in\mathcal F(B(R))\), Cegrell's integration by parts
formula \cite[Theorem~3.2]{Ceg04}, applied with \(u=g_{B(R)}(z,\cdot)\),
\(v=\varphi_z\), and \(T=(dd^c\varphi_z)^{n-1}\), gives
\[
\int_{B(R)} g_{B(R)}(z,w)\,(dd^c\varphi_z)^n(w)
\ge
\int_{B(R)} \varphi_z(w)\, dd^c g_{B(R)}(z,w)\wedge(dd^c\varphi_z)^{n-1}(w).
\tag{4.4}
\]
Since \(\varphi_z\ge\widetilde W_z\ge-\|u\|_\infty-\Lambda\) and
\(dd^c g_{B(R)}(z,\cdot)\wedge(dd^c\varphi_z)^{n-1}\ge0\), we get
\[
\int_{B(R)} g_{B(R)}(z,w)\,(dd^c\varphi_z)^n(w)
\ge
-(\|u\|_\infty+\Lambda)
\int_{B(R)} dd^c g_{B(R)}(z,w)\wedge(dd^c\varphi_z)^{n-1}(w).
\tag{4.5}
\]
By Corollary~5.6 of \cite{Ceg04}, applied with
\(u_1=\cdots=u_{n-1}=\varphi_z\) and \(u_n=g_{B(R)}(z,\cdot)\), and using
\((dd^c g_{B(R)}(z,\cdot))^n=(2\pi)^n\delta_z\) together with
\((dd^c\varphi_z)^n=\mathbf 1_{B(z,r_1)}\mu\), we obtain
\[\begin{array}{ll}\tag{4.6}\displaystyle
\int_{B(R)} dd^c g_{B(R)}(z,\cdot)\wedge(dd^c\varphi_z)^{n-1}
&\le
\displaystyle\left(\int_{B(R)}(dd^c g_{B(R)}(z,\cdot))^n\right)^{1/n}
\left(\int_{B(R)}(dd^c\varphi_z)^n\right)^{(n-1)/n}\\
&=
2\pi\,M^{(n-1)/n},
\end{array}
\]
where \(M:=\mu(\Omega)\).
Combining (4.5) and (4.6),
\[
\int_{B(R)} g_{B(R)}(z,w)\,(dd^c\varphi_z)^n(w)
\ge
-2\pi\,(\|u\|_\infty+\Lambda)\,M^{(n-1)/n}.
\tag{4.7}
\]
Since \((dd^c\varphi_z)^n=\mathbf 1_{B(z,r_1)}\mu\),
\[
\int_{B(z,r_1)} g_{B(R)}(z,w)\,d\mu(w)
\ge
-2\pi\,(\|u\|_\infty+\Lambda)\,M^{(n-1)/n}.
\tag{4.8}
\]
For the outer part, use
\(g_\Omega(z,w)\ge g_{B(R)}(z,w)\ge\log(|z-w|/(2R))\): for
\(w\in\Omega\setminus B(z,r_1)\), \(|z-w|\ge r_1\), hence
\[
\int_{\Omega\setminus B(z,r_1)}g_\Omega(z,w)\,d\mu(w)
\ge M\log\frac{r_1}{2R}.
\tag{4.9}
\]
Adding (4.8) and (4.9),
\[
p(z)\ge-2\pi\,(\|u\|_\infty+\Lambda)\,M^{(n-1)/n}
+M\log\frac{r_1}{2R}=:-C_K,
\tag{4.10}
\]
where \(C_K\) depends only on \(K\), \(\Omega\), \(R\), \(r_1\), \(r_2\),
\(\|u\|_\infty\), and \(M\), but not on \(z\in K\).

Since \(\Omega\) is hyperconvex,
then \(p_\mu\) vanishes on the boundary. There
exists \(\delta>0\) such that
\[
p(z)\ge-1\qquad\text{whenever }\operatorname{dist}(z,\partial\Omega)<\delta.
\]
Set \(K_\delta:=\{z\in\Omega:\operatorname{dist}(z,\partial\Omega)\ge\delta\}\), so
\(K_\delta\Subset\Omega\). Applying (4.10) with \(K=K_\delta\) gives
\(p\ge-C_{K_\delta}\) on \(K_\delta\), while \(p\ge-1\) on
\(\Omega\setminus K_\delta\). Hence
\[
p(z)\ge-\max\{C_{K_\delta},1\}=:-C
\qquad\text{for all }z\in\Omega.
\]
Together with \(p\le0\), this proves that \(p\) is bounded.

\medskip\noindent\textbf{Step 2: uniform control of the Monge--Amp\`ere mass.}
By Lemma~\ref{lem:representation} applied to the finite measure \(\mu\),
\[
(dd^c p)^n
=
\int_{\Omega^n}
dd^c_z g_\Omega(z,v_1)\wedge\cdots\wedge dd^c_z g_\Omega(z,v_n)\,
d\mu(v_1)\cdots d\mu(v_n).
\]
For every \((v_1,\dots,v_n)\in\Omega^n\), the functions
\(g_\Omega(\cdot,v_i)\) belong to \(\mathcal F(\Omega)\)  and
\((dd^c g_\Omega(\cdot,v_i))^n=(2\pi)^n\delta_{v_i}\). By Corollary~5.6 of
\cite{Ceg04}, applied to \(u_i=g_\Omega(\cdot,v_i)\), we have
\[
\int_\Omega dd^c g_\Omega(\cdot,v_1)\wedge\cdots\wedge dd^c g_\Omega(\cdot,v_n)
\le
\prod_{i=1}^n
\left(\int_\Omega(dd^c g_\Omega(\cdot,v_i))^n\right)^{1/n}
=
(2\pi)^n.
\tag{4.11}
\]
Combining the representation formula with (4.11) and Fubini,
\[
\int_\Omega(dd^c p)^n
\le (2\pi)^n\int_{\Omega^n}d\mu(v_1)\cdots d\mu(v_n)
=(2\pi)^n\,\mu(\Omega)^n<+\infty .
\tag{4.12}
\]
Combining Steps 1 and 2 with \(p\in\operatorname{PSH}(\Omega)\), we conclude that
\(p\in L^\infty(\Omega)\), \(p\) vanishes on \(\partial\Omega\), and
\(\int_\Omega(dd^c p)^n<+\infty\). Therefore \(p\in\mathcal E_0(\Omega)\). This
proves Theorem~\ref{thm:stability}(a).

\begin{remark}
\label{rem:carlehed}
Carlehed \cite[Proposition~9.3]{Car99} proved the \(\mathcal E_0\) statement
under the additional assumption that the Green function of \(\Omega\) is
\emph{symmetric}, i.e.\ \(g_\Omega(z,w)=g_\Omega(w,z)\). This assumption is
satisfied for convex domains by Lempert's theorem \cite{Lem81}, but it fails
in general. Bedford and Demailly \cite[Proposition~2]{BD88} constructed a
bounded strictly pseudoconvex domain \(\Omega\subset\mathbb C^2\) with real-analytic
boundary whose Green function is not symmetric. Their example is
\[
\Omega=\Bigl\{(z,w)\in\mathbb C^2:\ 
\max\Bigl(\tfrac12\log\frac{|w^2-z^2(z-a)|}{\varepsilon^2},\ \log|z|\Bigr)<0
\Bigr\},
\]
with \(0<|a|<1\) and \(\varepsilon>0\) small, for which
\(g_\Omega((0,0),(a,0))\ne g_\Omega((a,0),(0,0))\). The symmetry is used in
\cite{Car99} to make \(g_\Omega(z,\cdot)\) plurisubharmonic in the second
variable, so that the test functions \(\max\{g_\Omega(z,\cdot),N\}\) can be
used in the energy estimates. Our proof avoids this issue by replacing
\(g_\Omega(z,\cdot)\) by the explicit Green function \(g_{B(R)}(z,\cdot)\) of a
ball \(B(R)\Supset\Omega\) in the integration by parts: the comparison
\(g_\Omega\ge g_{B(R)}\) provides the required lower bound, and \(g_{B(R)}\) is
symmetric, hence amenable to the standard integration by parts. Therefore
Theorem~\ref{thm:stability}(a) extends Carlehed's Proposition~9.3 to arbitrary
bounded strictly pseudoconvex domains with \(C^\infty\) boundary, without any
symmetry assumption on the Green function.
\end{remark}

\subsection{Proof of Theorem~\ref{thm:stability}(b): the \(\mathcal F_1\) case}

Set \(E:=\int_\Omega(-u)(dd^c u)^n<\infty\) and \(M:=\mu(\Omega)<\infty\).
By Cegrell's theorem \cite[Theorem~1.1]{Ceg12}, \(p_\mu\in\mathcal F(\Omega)\).
By definition of \(\mathcal F(\Omega)\), there exists a decreasing sequence
\((p_j)\subset\mathcal E_0(\Omega)\) such that \(p_j\searrow p_\mu\) on \(\Omega\) and
\[
\sup_j\int_\Omega(dd^c p_j)^n<+\infty.
\tag{4.13}
\]
By the sharp estimate \cite[Theorem~1.1]{Ceg12} applied to the finite
measure \(\mu\) and \(h=p_j\in\operatorname{PSH}^-(\Omega)\) and Hölder's inequality
(with \(M=\mu(\Omega)\)) we get
\[
\begin{array}{ll}
     A_j:=\displaystyle \int_\Omega(-p_j)(dd^c p_j)^n
&\le
\displaystyle\int_\Omega(-p_j)(dd^c p_\mu)^n \tag{4.14}\\ &\le \displaystyle\left(\int_\Omega(-p_j)^{1/n}\,d\mu\right)^n\\ &\le \displaystyle M^{n-1}\int_\Omega(-p_j)\,d\mu.
\end{array}
\]
Since \(u\in\mathcal F_1(\Omega)\subset\mathcal E_1(\Omega)\), there exists a decreasing sequence
\((u_k)\subset\mathcal E_0(\Omega)\) with \(u_k\searrow u\). By Theorem~3.8 of
\cite{Ceg98},
\[
E_k:=\int_\Omega(-u_k)(dd^c u_k)^n\nearrow E:=\int_\Omega(-u)(dd^c u)^n,
\tag{4.15}
\]
so \(E_k\le E\) for all \(k\).

For each \(k\), since \(p_j\in\mathcal E_0(\Omega)\subset\mathcal E_1(\Omega)\) and
\(u_k\in\mathcal E_0(\Omega)\), Theorem~3.2 of \cite{Ceg98} gives
\[
\begin{array}{ll}
\displaystyle
\int_\Omega(-p_j)(dd^c u_k)^n
&\le \displaystyle
\left(\int_\Omega(-p_j)(dd^c p_j)^n\right)^{1/(n+1)}
\left(\int_\Omega(-u_k)(dd^c u_k)^n\right)^{n/(n+1)}
\\
&=\displaystyle
A_j^{1/(n+1)}E_k^{n/(n+1)}\\
&\le
A_j^{1/(n+1)}E^{n/(n+1)}.
\end{array}
\tag{4.16}
\]
By Proposition~5.1 of \cite{Ceg04}, since \(u_k\searrow u\) with
\(u_k\in\mathcal E_0\subset\mathcal F\) and \(u\in\mathcal F_1\subset\mathcal F\),
and since \(-p_j\in\operatorname{PSH}^-(\Omega)\), the Monge--Amp\`ere measures
\((dd^c u_k)^n\) satisfy
\(
\lim_{k\to\infty}\int_\Omega(-p_j)(dd^c u_k)^n
=
\int_\Omega(-p_j)(dd^c u)^n.
\)
Hence
\[
\int_\Omega(-p_j)\,d\mu=\int_\Omega(-p_j)(dd^c u)^n
\le
A_j^{1/(n+1)}E^{n/(n+1)}.
\tag{4.17}
\]
Combining (4.14) and (4.17), we get 
\(
A_j\le M^{n-1}A_j^{1/(n+1)}E^{n/(n+1)}.
\)
If \(A_j>0\), dividing by \(A_j^{1/(n+1)}\) yields
\(
A_j^{n/(n+1)}\le M^{n-1}E^{n/(n+1)}
\), hence
\[
A_j\le M^{(n-1)(n+1)/n}E=:C<\infty.
\tag{4.18}
\]
We have \(p\in\mathcal F(\Omega)\)  and a sequence \(p_j\in\mathcal E_0(\Omega)\)
with \(p_j\searrow p\) and \(\sup_j\int(-p_j)(dd^c p_j)^n\le C<\infty\). By the definition of \(\mathcal E_1(\Omega)\) \cite[Definition~2.2]{Ceg98},
this gives \(p\in\mathcal E_1(\Omega)\). Since \(\mathcal F_1=\mathcal F\cap\mathcal E_1\), we conclude
\(p\in\mathcal F_1(\Omega)\). This proves Theorem~\ref{thm:stability}(b).

\begin{remark}

\label{rem:R1}
In \cite{Ceg12}, Cegrell obtains \(p_\nu\in\mathcal F_1(\Omega)=\mathcal{F}(\Omega)\cap \mathcal{E}_1(\Omega)\) under the domination
condition \[0<\nu\le(-u)^{(n-1)/n}(dd^c v)^n\qquad \mbox{with}\; u,v\in\mathcal E_1(\Omega).\]
Our Theorem~\ref{thm:stability}(b) treats a different class of source
measures, namely \(\mu=(dd^c u)^n\) with \(u\in\mathcal F_1(\Omega)\). 
Thus, the two results provide sufficient conditions of
different types for the Green potential to belong to
\(\mathcal F_1(\Omega)\). We do not establish any implication between
these two conditions, nor do we claim that they are
incomparable.

\end{remark}

\section{Infinite-mass sources}
\label{sec:infinite-mass}

In this section we consider sources \(\mu\) of infinite total mass.  In
Section~\ref{sec:infinite-mass-example} we construct an explicit example
where \(\mu(\mathbb B)=+\infty\) and \(p_\mu\in\mathcal F(\mathbb B)\), showing that the smoothing
phenomenon is not restricted to finite sources.  In
Section~\ref{sec:infinite-mass-general} we prove that, for every
\(u\in\mathcal E_1(\mathbb B)\) with \(\mu=(dd^c u)^n\) possibly of infinite mass, the
Green potential \(p_\mu\) belongs to \(\mathcal N(\mathbb B)\), and we formulate two open
questions about the energy class \(\mathcal E_1\) and the absolute continuity of
\((dd^c p_\mu)^n\).

\subsection{An explicit example}
\label{sec:infinite-mass-example}

We construct \(\mu\) on \(\mathbb B\subset\mathbb C^2\) with \(\mu(\mathbb B)=+\infty\), whose
Green potential \(u=p_\mu\) is finite, locally bounded, belongs to
\(\mathcal F(\mathbb B)\), and satisfies
\[
(dd^c u)^2\ll\lambda,
\qquad
\int_\mathbb B(dd^c u)^2<+\infty.
\]
For \(k\ge2\), set \(r_k=1-2^{-k}\), \(S_{r_k}=\partial B(0,r_k)\), and let
\(\sigma_{r_k}\) be the Euclidean surface measure on \(S_{r_k}\). Define
\[
\mu:=\sum_{k=2}^\infty k\,\sigma_{r_k}.
\tag{5.1}
\]

\begin{proposition}
\label{prop:mu-radon}
The measure \(\mu\) is a positive atomless Radon measure on \(\mathbb B\),
\[
\int_\mathbb B(1-|w|)\,d\mu(w)<+\infty,
\qquad
\mu(\mathbb B)=+\infty.
\]
\end{proposition}

\begin{proof}
Local finiteness follows from \(r_k\to1\); atomlessness from that of each
\(\sigma_{r_k}\). Since \(\sigma_{r_k}(S_{r_k})=2\pi^2r_k^3\), we have
\(\mu(\mathbb B)=2\pi^2\sum_{k\ge2}kr_k^3=+\infty\), while \(1-r_k=2^{-k}\) gives
\[
\int_\mathbb B(1-|w|)\,d\mu(w)
\le 2\pi^2\sum_{k\ge2}k\,2^{-k}<+\infty.
\]
\end{proof}

Carlehed's formula \cite[Proposition~3.2]{Car99} gives, with \(t=|z|\),
\[
p_{\sigma_r}(z)=\pi^2r^3
\begin{cases}
\displaystyle
\log r^2+\frac{1-r^2}{r^2}
+\frac{(1-r^2)(1-t^2)}{r^2t^2}\log(1-t^2),
&0\le t\le r,\\[1.1em]
\displaystyle
\log t^2+\frac{1-t^2}{t^2}
+\frac{(1-r^2)(1-t^2)}{r^2t^2}\log(1-r^2),
&r\le t<1.
\end{cases}
\tag{5.2}
\]
In particular, \(p_{\sigma_r}\) is continuous on \(\overline{\mathbb B}\), bounded,
and vanishes on \(\partial\mathbb B\).

For \(N\ge2\), put
\[
u_N:=\sum_{k=2}^N k\,p_{\sigma_{r_k}},
\qquad
u:=\sum_{k=2}^\infty k\,p_{\sigma_{r_k}}.
\]

\begin{lemma}
\label{lem:decay}
For every \(R\in(0,1)\), there is \(C_R>0\) such that, for \(r>R\) close to
\(1\),
\[
\sup_{|z|\le R}|p_{\sigma_r}(z)|\le C_R(1-r^2).
\]
\end{lemma}

\begin{proof}
For \(t\le R\) and \(s=1-r^2\), the first branch of (5.2) reads
\[
\log(1-s)+\frac{s}{1-s}+\frac{s}{1-s}h(t),
\qquad
h(t)=t^{-2}(1-t^2)\log(1-t^2),
\]
with \(h\) bounded on \([0,R]\) and \(h(0)=-1\). Since
\(\log(1-s)+s/(1-s)=O(s^2)\) and \(s/(1-s)=O(s)\), the estimate follows.
\end{proof}

\begin{proposition}
\label{prop:u-psh}
The series defining \(u\) converges uniformly on compact subsets of \(\mathbb B\),
and \(u=p_\mu\). In particular, \(u\) is finite, negative, and
plurisubharmonic on \(\mathbb B\).
\end{proposition}

\begin{proof}
On \(K\subset\{|z|\le R\}\Subset\mathbb B\), Lemma~\ref{lem:decay} gives
\(\sup_K k|p_{\sigma_{r_k}}|\le C_R k(1-r_k^2)\le C_R' k\,2^{-k}\), which is
summable. Hence \(u\) is a locally uniform limit of psh functions, so
\(u\in\operatorname{PSH}(\mathbb B)\), \(u\le0\), and \(u\) is finite. Tonelli applied to
\(-g_\mathbb B(z,\cdot)\ge0\) yields
\(
\int_\mathbb B g_\mathbb B(z,w)\,d\mu(w)
=\sum_{k\ge2}k\,p_{\sigma_{r_k}}(z)=u(z).
\)
Thus \(u=p_\mu\).
\end{proof}

Write \(p_k:=p_{\sigma_{r_k}}\).

\begin{lemma}
\label{lem:mixed-representation}
For \(k,l\ge2\),
\[
dd^c p_k\wedge dd^c p_l
=
\int_{S_{r_k}\times S_{r_l}}
dd^c_z g_\mathbb B(z,v_1)\wedge
dd^c_z g_\mathbb B(z,v_2)\,
d\sigma_{r_k}(v_1)d\sigma_{r_l}(v_2).
\tag{5.3}
\]
\end{lemma}

\begin{proof}
Apply the representation formula (3.1) in dimension two to
\(\sigma_{r_k}+\sigma_{r_l}\), expand, and subtract the two identities
corresponding to \(\sigma_{r_k}\) and \(\sigma_{r_l}\).
\end{proof}

\begin{proposition}
\label{prop:u-AC-no-lemma}
In the sense of positive measures on \(\mathbb B\),
\[
(dd^c u)^2
=
\sum_{k,l\ge2}kl\,dd^c p_k\wedge dd^c p_l.
\tag{5.4}
\]
Every \(dd^c p_k\wedge dd^c p_l\) is absolutely continuous with respect to
\(\lambda\), the series converges in total mass, and consequently
\[
(dd^c u)^2\ll\lambda,
\qquad
\int_\mathbb B(dd^c u)^2<+\infty.
\tag{5.5}
\]
\end{proposition}

\begin{proof}
Since \(u_N\searrow u\) locally uniformly and \(u\) is locally bounded, we
have \((dd^c u_N)^2\longrightarrow(dd^c u)^2\) weakly on \(\mathbb B\). By
multilinearity,
\((dd^c u_N)^2=\sum_{k,l=2}^N kl\,dd^c p_k\wedge dd^c p_l\). Set
\[
G_{k,l}:=\int_\mathbb B dd^c p_k\wedge dd^c p_l.
\]
Since \(g_\mathbb B(\cdot,v_j)\in \mathcal F(\mathbb B)\), for \(j=1,2\), and the Green function of the ball is symmetric
(\(g_\mathbb B(v_2,v_1)=g_\mathbb B(v_1,v_2)\)), then
\[
\int_\mathbb B dd^c_z g_\mathbb B(z,v_1)\wedge dd^c_z g_\mathbb B(z,v_2)
=-(2\pi)^2g_\mathbb B(v_1,v_2),
\qquad v_1\ne v_2.
\]
Together with (5.3) and the fact that the diagonal has zero product surface
measure, we obtain
\(G_{k,l}=-(2\pi)^2\int_{S_{r_k}\times S_{r_l}}g_\mathbb B(v_1,v_2)\,d\sigma_{r_k}(v_1)d\sigma_{r_l}(v_2)\).
For \(k\le l\), integrating first in \(v_1\) gives
\(G_{k,l}=-(2\pi)^2\,2\pi^2r_l^3\,p_k(r_l)\). A direct computation from
(5.2), with \(s_j=1-r_j^2\asymp2^{-j}\), yields
\[
0\le G_{k,l}\le C\min(k,l)\,2^{-k-l},
\qquad k,l\ge2.
\tag{5.6}
\]
Therefore
\(\sum_{k,l\ge2}kl\,G_{k,l}\le C\sum_{k,l\ge2}kl\min(k,l)\,2^{-k-l}<+\infty\),
so \(\nu:=\sum_{k,l\ge2}kl\,dd^c p_k\wedge dd^c p_l\) converges in total
mass. For \(\chi\in C_c^\infty(\mathbb B)\), this and multilinearity give
\(\int_\mathbb B\chi\,(dd^c u_N)^2\to\sum_{k,l\ge2}kl\int_\mathbb B\chi\,dd^c p_k\wedge dd^c p_l\).
Comparing with weak convergence proves (5.4), and the summability gives
\(\int_\mathbb B(dd^c u)^2<+\infty\).

It remains to prove absolute continuity. If \(k\ne l\), the spheres are
disjoint, so \(v_1\ne v_2\) on \(S_{r_k}\times S_{r_l}\). By
Proposition~\ref{prop:mixed},
\(dd^c_z g_\mathbb B(z,v_1)\wedge dd^c_z g_\mathbb B(z,v_2)\ll\lambda\), hence
\(dd^c p_k\wedge dd^c p_l\ll\lambda\). If \(k=l\), the diagonal
\(\Delta_k=\{(v,v):v\in S_{r_k}\}\) is
\(\sigma_{r_k}\otimes\sigma_{r_k}\)-null, and the same argument applies off
\(\Delta_k\), so \((dd^c p_k)^2\ll\lambda\). Since the series converges in
total mass, its sum is absolutely continuous. This proves (5.5).
\end{proof}

\begin{proof}[Completion of the example]
For every \(N\), \(u_N\in\mathcal E_0(\mathbb B)\): it is bounded, psh, vanishes on
\(\partial\mathbb B\), and by (5.6),
\[
\int_\mathbb B(dd^c u_N)^2
=\sum_{k,l=2}^N kl\,G_{k,l}
\le C\sum_{k,l\ge2}kl\min(k,l)\,2^{-k-l}<+\infty,
\]
uniformly in \(N\). Since \(u_N\searrow u\), the definition of \(\mathcal F(\mathbb B)\)
gives \(u\in\mathcal F(\mathbb B)\).
\end{proof}

\begin{theorem}
\label{thm:main-infinite}
For \(\mu\) defined by (5.1), its Green potential \(u=p_\mu\) satisfies
\[
u\in\mathcal F(\mathbb B),
\qquad
(dd^c u)^2\ll\lambda,
\qquad
\int_\mathbb B(dd^c u)^2<+\infty,
\]
whereas \(\mu(\mathbb B)=+\infty\).
\end{theorem}

\begin{remark}
A positive Radon measure of infinite total mass may thus have a finite,
locally bounded Green potential in \(\mathcal F(\mathbb B)\), with finite and absolutely
continuous Monge--Amp\`ere measure, showing that the smoothing phenomenon is
not restricted to finite sources.
\end{remark}

\subsection{General infinite-mass sources}
\label{sec:infinite-mass-general}

We now turn to the general case: \(u\in\mathcal E_1(\mathbb B)\),
\(\mu:=(dd^c u)^n\) possibly of infinite mass, and
\(p_\mu:=p_\mu\). The following pointwise estimate does not
involve the total mass \(\mu(\mathbb B)\).

\begin{theorem}

Let \(u\in\mathcal E_1(\mathbb B)\), \(\mu=(dd^c u)^n\), and
\(p_\mu=p_\mu\). Then there exists a constant \(C=C(n)>0\)
such that
\(
-p_\mu(z)\le
C E^{(n-1)/n}(-u(z))^{1/n}\),
\( z\in\mathbb B,
\)
where
\(
E:=\int_\mathbb B(-u)(dd^c u)^n<+\infty.
\)
At points where \(u(z)=-\infty\), the right-hand side is
understood to be \(+\infty\). Moreover \(p_\mu\in\mathcal{N}(\mathbb B).\)
\end{theorem}

\begin{proof}
Choose a decreasing sequence \(u_k\in\mathcal E_0(\mathbb B)\) such that
\(
u_k\searrow u.
\)
Set
\(
\mu_k:=(dd^c u_k)^n\), \(p_k:=p_{\mu_k}\) and \(E_k:=\int_\mathbb B(-u_k)(dd^c u_k)^n.
\)
By the standard energy approximation, we have
\(
E_k\le E.
\)
We first establish the monotonicity of the sequence
\((p_k)\). Fix \(z\in\mathbb B\). Since the pluricomplex Green
function of the unit ball is symmetric,
\(
g_\mathbb B(z,w)=g_\mathbb B(w,z),
\)
the function \(w\mapsto g_\mathbb B(z,w)\) is negative
plurisubharmonic in \(\mathbb B\). 
Let \(j\ge k\). Since \(u_j\le u_k\), then 

\[ -p_k(z) =\displaystyle\int_\Omega -g(z,w)(dd^c u_k)^n \leq \displaystyle\int_\Omega -g(z,w)(dd^c u_j)^n=-p_j(z).\]
Thus \((p_k)\) is decreasing. We next apply integration by parts in the Cegrell classes.
Since \(g_\mathbb B(z,\cdot)\in\mathcal F(\mathbb B)\) and
\(
(dd^c g_\mathbb B(z,\cdot))^n=(2\pi)^n\delta_z,
\)
we obtain
\[
-p_k(z)
=\int_\mathbb B(-g_\mathbb B(z,w))\,d\mu_k(w)
=\int_\mathbb B(-u_k)\,
  dd^c g_\mathbb B(z,\cdot)\wedge(dd^c u_k)^{n-1}.
\]
Cegrell's H\"older inequality, applied with
\(u_1=\cdots=u_{n-1}=u_k\) and
\(u_n=g_\mathbb B(z,\cdot)\), yields
\[
-p_k(z)
\le
\left(
\int_\mathbb B(-u_k)(dd^c g_\mathbb B(z,\cdot))^n
\right)^{1/n}
\times
\left(
\int_\mathbb B(-u_k)(dd^c u_k)^n
\right)^{(n-1)/n}.
\]
By the normalization of the Green function, the first
factor is
\(
2\pi(-u_k(z))^{1/n},
\)
whereas the second factor is \(E_k^{(n-1)/n}\).
Therefore,
\[
-p_k(z)
\le
2\pi E^{(n-1)/n}(-u_k(z))^{1/n}.
\tag{5.7}
\]
It remains to identify the limit of \(p_k\). By Proposition~5.1 of \cite{Ceg04}, applied to the decreasing sequence
\((u_k)\) and the negative plurisubharmonic test function
\(g_\mathbb B(z,\cdot)\), we have
\[
\begin{aligned}
\lim_{k\to\infty}
\int_\mathbb B(-g_\mathbb B(z,w))(dd^c u_k)^n
=
\int_\mathbb B(-g_\mathbb B(z,w))(dd^c u)^n.
\end{aligned}
\]
It follows that
\(
p_k(z)\searrow p_\mu(z),
\)
where the limit is understood in the extended sense and
may equal \(-\infty\).
Hence by  (5.7),  we obtain
\(
-p_\mu(z)
\le
2\pi E^{(n-1)/n}(-u(z))^{1/n}.
\)

Now, since \(u\in\mathcal E_1(\mathbb B)\subset\mathcal N(\mathbb B)\), and
\(\mathcal N(\mathbb B)\) is stable under the transformation
\(
u\longmapsto -C(-u)^\alpha\) for \(0<\alpha<1\), (this stability follows directly from the definition
of \(\mathcal N(\mathbb B)\)),
we have
\[
v:=-C E^{(n-1)/n}(-u)^{1/n}\in\mathcal N(\mathbb B).
\]
Since \(
v\le p_\mu
\), then \(\widetilde{p_\mu}=0\). Therefore,
\(
{p_\mu\in\mathcal N(\mathbb B).}
\)

\end{proof}

\begin{remark}
\label{rem:open}

Two natural questions remain open in the infinite-mass setting.

\medskip
\noindent\emph{1. The energy class \(\mathcal E_1\).}
We conjecture that \(p_\mu\in\mathcal E_1(\mathbb B)\), namely,
\[
\int_\mathbb B(-p_\mu)(dd^c p_\mu)^n<+\infty.
\tag{5.8}
\]
While this holds in the finite-mass case by
Theorem~\ref{thm:stability}(b), the standard energy estimates
in the infinite-mass case involve the potentially unbounded
factor \(\mu_j(\mathbb B)^{n-1}\).

\medskip
\noindent\emph{2. Absolute continuity of the Monge--Amp\`ere
measure.}
We also conjecture that
\[
(dd^c p_\mu)^n\ll\lambda.
\tag{5.9}
\]
This is established in the finite-mass case by
Theorem~\ref{thm:bdd}. In the infinite-mass case, the
monotonicity \(p_k\searrow p_\mu\), established above,
and the continuity of the Monge--Amp\`ere operator  yield
\[
(dd^c p_k)^n\rightharpoonup(dd^c p_\mu)^n.
\tag{5.10}
\]
However, weak convergence does not preserve absolute
continuity. Thus, proving (5.9) requires an additional
argument, such as establishing uniform absolute continuity
of the approximating Monge--Amp\`ere measures.

\end{remark}


\begin{thebibliography}{99}

\bibitem{ACC09}
P.~Åhag, U.~Cegrell, R.~Czy\.z and P.~H.~Hi\^ep,
\emph{Monge--Amp\`ere measures on pluripolar sets},
J. Math. Pures Appl. \textbf{92} (2009), 613--627.

\bibitem{BT76}
E.~Bedford and B.~A. Taylor,
\emph{The Dirichlet problem for a complex Monge--Amp\`ere equation},
Invent. Math. \textbf{37} (1976), 1--44.

\bibitem{BD88}
E.~Bedford and J.-P. Demailly,
\emph{Two counterexamples concerning the pluricomplex Green function in
\(\mathbb C^n\)},
Indiana Univ. Math. J. \textbf{37} (1988), 865--867.

\bibitem{Blo00}
Z.~B\l ocki,
\emph{The \(C^{1,1}\) regularity of the pluricomplex Green function},
Michigan Math. J. \textbf{47} (2000), 211--214.

\bibitem{Blo04}
Z.~B\l ocki,
\emph{On the definition of the Monge--Amp\`ere operator in \(\mathbb C^2\)},
Math. Ann. \textbf{328} (2004), 415--423.

\bibitem{Blo06}
Z.~B\l ocki,
\emph{The domain of definition of the complex Monge--Amp\`ere operator},
Amer. J. Math. \textbf{128} (2006), 519--530.

\bibitem{Car99}
M.~Carlehed,
\emph{Potentials in pluripotential theory},
Ann. Fac. Sci. Toulouse Math. (6) \textbf{8} (1999), no.~3, 439--469.

\bibitem{Ceg98}
U.~Cegrell,
\emph{Pluricomplex energy},
Acta Math. \textbf{180} (1998), 187--217.

\bibitem{Ceg04}
U.~Cegrell,
\emph{The general definition of the complex Monge--Amp\`ere operator},
Ann. Inst. Fourier (Grenoble) \textbf{54} (2004), no.~1, 159--179.

\bibitem{Ceg12}
U.~Cegrell,
\emph{Measures of finite pluricomplex energy},
Ann. Polon. Math. \textbf{106} (2012), 107--112.

\bibitem{CzyzHed}
R.~Czy\.z and L.~Hed,
\emph{Subextension of plurisubharmonic functions without increasing the
total Monge--Amp\`ere mass},
Ann. Polon. Math. \textbf{94} (2008), 275--281.

\bibitem{Dem87}
J.-P. Demailly,
\emph{Mesures de Monge--Amp\`ere et mesures pluriharmoniques},
Math. Z. \textbf{194} (1987), 519--564.

\bibitem{Dem93}
J.-P. Demailly,
\emph{Monge--Amp\`ere operators, Lelong numbers and intersection theory},
Complex Analysis and Geometry, Univ. Series in Math., edited by V.~Ancona
and A.~Silva, Plenum Press, New York (1993).

\bibitem{GZ17}
V.~Guedj and A.~Zeriahi,
\emph{Degenerate Complex Monge--Amp\`ere Equations},
EMS Tracts in Mathematics 26, European Mathematical Society, 2017.

\bibitem{Hed}
L.~Hed,
\emph{Approximation of negative plurisubharmonic functions with given
boundary values},
Internat. J. Math. \textbf{21} (2010), no.~9, 1135--1145.

\bibitem{Kli91}
M.~Klimek,
\emph{Pluripotential Theory},
London Math. Soc. Monogr. (N.S.) 6, Clarendon Press, Oxford, 1991.

\bibitem{Lem81}
L.~Lempert,
\emph{La m\'etrique de Kobayashi et la repr\'esentation des domaines sur la
boule},
Bull. Soc. Math. France \textbf{109} (1981), 427--474.

\bibitem{KHH09}
N.~V. Khue and P.~H. Hiep,
\emph{A comparison principle for the complex Monge--Amp\`ere operator in
Cegrell's classes and applications},
Proc. Amer. Math. Soc. \textbf{137} (2009), 553--564.

\bibitem{Wik05}
J.~Wiklund,
\emph{Pluricomplex charge at weak singularities},
preprint arXiv:math/0510671 (2005).

\end{thebibliography}
\end{document}